# MULTI-VALUED BOUNDARY VALUE PROBLEMS INVOLVING LERAY-LIONS OPERATORS AND DISCONTINUOUS NONLINEARITIES

SIMONA DĂBULEANU – VICENŢIU RĂDULESCU

We prove an existence result for a class of Dirichlet boundary value problems with discontinuous nonlinearity and involving a Leray-Lions operator. The proof combines monotonicity methods for elliptic problems, variational inequality techniques and basic tools related to monotone operators. Our work generalizes a result obtained in Carl [4].

## 1. Introduction and the main result.

Let $\Omega \subset \mathbf{R}^N$ be a bounded domain with smooth boundary. Consider the boundary value problem

$$(P)\begin{cases} -\mathrm{div}\,(a(x,\nabla u(x))) = f(u(x)), & \text{if } x \in \Omega \\ u = 0, & \text{on } \partial\,\Omega, \end{cases}$$

where $a : \Omega \times \mathbf{R}^N \to \mathbf{R}^N$ is a Carathéodory function having the properties

$(a_1)$ there exist $p > 1$ and $\lambda > 0$ such that $a(x,\xi)\cdot\xi \geq \lambda \cdot \|\xi\|^p$, for a.e. $x \in \Omega$ and for any $\xi \in \mathbf{R}^N$;

$(a_2)$ $(a(x,\xi) - a(x,\eta)) \cdot (\xi - \eta) > 0$, for any $\xi, \eta \in \mathbf{R}^N$, $\xi \neq \eta$;

$(a_3)$ there exist $\alpha \in \mathbf{R}^+$ and $k \in L^{p'}(\Omega)$ such that $|a(x,\xi)| \leq \alpha(k(x)+|\xi|^{p-1})$, for a.e. $x \in \Omega$ and for any $\xi \in \mathbf{R}^N$.

Assume that the nonlinearity $f : \mathbf{R} \to \mathbf{R}$ satisfies the hypothesis

$(H_1)$ there exist nondecreasing functions $f, g : \mathbf{R} \to \mathbf{R}$ such that $f = g - h$ .

Let $\beta : \mathbf{R} \to 2^{\mathbf{R}}$ be the maximal monotone graph associated with the nondecreasing function $h$ (see Brezis [3]). More exactly,

$$\beta(s) := [h^-(s), h^+(s)], \qquad \text{for all } s \in \mathbf{R},$$

where

$$h^-(s) = \lim_{\varepsilon \to 0+} h(s-\varepsilon), \quad h^+(s) = \lim_{\varepsilon \to 0+} h(s+\varepsilon)\,.$$

Under this assumption we reformulate the problem $(P)$ as follows

$$(P') \begin{cases} -\mathrm{div}\,(a(x, \nabla(x))) + \beta(u(x)) \ni g(u(x)), & \text{if } x \in \Omega \\ u = 0, & \text{on } \partial\Omega. \end{cases}$$

Denote by $G$ the Nemitskii operator associated with $g$, that is, $G(u)(x) = g(u(x))$.

DEFINITION 1. *A function $u \in W_0^{1,p}(\Omega)$ is called a solution of the problem $(P')$ if there exists $v \in L^{p'}(\Omega)$ such that*

*i)* $v(x) \in \beta(u(x))$ *a.e. in* $\Omega$,

*ii)* $\int_\Omega a(x, \nabla u)\cdot \nabla w\, dx + \int_\Omega v\cdot w\, dx = \int_\Omega G(u)\cdot w\, dx$, *for any* $w \in W_0^{1,p}(\Omega)$.

Let $L_+^p$ be the set of nonnegative elements of $L^p(\Omega)$. For any $v, w \in \Omega$ such that $v \le w$, we set

$$[v, w] = \{u \in L^p(\Omega) \;/\; v \le u \le w\}.$$

DEFINITION 2. *A function $\overline{u} \in W^{1,p}(\Omega)$ is called an upper solution of the problem $(P')$ if there exists a function $\overline{v} \in L^{p'}(\Omega)$ such that*

*i)* $\overline{v}(x) \in \beta(\overline{u}(x))$ *a.e. in* $\Omega$,

*ii)* $\overline{u} \ge 0$ *on* $\partial\Omega$,

*iii)* $\int_\Omega a(x, \nabla\overline{u}) \cdot \nabla w\, dx + \int_\Omega \overline{v} \cdot w\, dx \ge \int_\Omega G(\overline{u}) \cdot w\, dx$ *for all* $w \in W_0^{1,p}(\Omega) \cap L_+^p(\Omega)$.

DEFINITION 3. *A function $\overline{u} \in W^{1,p}(\Omega)$ is called a lower solution of the problem $(P')$ if there exists a function $\overline{v} \in L^{p'}(\Omega)$ such that*

*i)* $\underline{v}(x) \in \beta(\underline{u}(x))$ *a.e. in* $\Omega$,

*ii)* $\underline{u} \le 0$ *on* $\partial\Omega$,

*iii)* $\int_\Omega a(x, \nabla \underline{u}) \cdot \nabla w\, dx + \int_\Omega \underline{v} \cdot w\, dx \le \int_\Omega G(\underline{u}) \cdot w\, dx$ *for any* $w \in W_0^{1,p}(\Omega) \cap L_+^p(\Omega)$.

In the sequel the following hypothesis will be needed:

$(H_2)$ There exist an upper solution $\overline{u}$ and a lower solution $\underline{u}$ of the problem $(P')$ such that $\underline{u} \le \overline{u}$, and $G(\underline{u}),\ G(\overline{u}),\ H^+(\overline{u}),\ H^-(\underline{u}) \in L^{p'}(\Omega)$.

The following is a generalization of the main result in Carl [4].

THEOREM 1. *Assume hypothesis* $(H_1)$ *and* $(H_2)$ *hold and that* $g$ *is right (resp. left) continuous. Then there exists a maximal (resp. minimal) solution* $u \in [\underline{u}, \overline{u}]$ *of the problem* $(P')$.

## 2. Proof of Theorem 1.

We first reformulate the problem (P') in terms of variational inequalities using the subdifferential theory in the sense of convex analysis.

Let $j : \mathbf{R} \to (-\infty, \infty]$ be a convex, proper and lower semicontinuous function. Let $\partial j$ be the subdifferential of $j$, that is

(1) $$\partial j(r) = \{\hat{r} \in \mathbf{R} : j(s) \ge j(r) + \hat{r}(s-r) \qquad \text{for all } s \in \mathbf{R}\}.$$

We recall the following result concerning maximal monotone graphs in $\mathbf{R}^2$ (see Brezis [3] [Corollary 2.10], p. 43)

LEMMA 1. *Let* $\beta : \mathbf{R} \to 2^{\mathbf{R}}$ *be a maximal monotone graph in* $\mathbf{R}^2$. *Then there exists a convex, proper and lower semicontinuous function* $j : \mathbf{R} \to (-\infty, +\infty]$ *such that* $\beta = \partial j$. *Moreover, the function* $j$ *is uniquely determined up to an additive constant.*

We observe that the function $h$ appearing in $(H_1)$ can always be chosen so that $h(0) = 0$. Then the maximal monotone graph $\beta$ has the properties

(2) $$D(\beta) = \mathbf{R} \quad \text{and} \quad 0 \in \beta(0)\,.$$

Since the function $j$ related to $\beta$ according to Lemma 1 is uniquely determined up to an additive constant we can assume that

(3) $$j(0) = 0\,.$$

So, by (1) , (2) and (3) it follows that

$$j(s) \geq 0 \qquad \text{for all } s \in \mathbf{R}. \tag{4}$$

Define $J : L^p(\Omega) \to (-\infty, +\infty]$ by

$$J(v) = \begin{cases} \int_\Omega j(v(x))\,dx, & \text{for } j(v(\cdot)) \in L^1(\Omega) \\ +\infty & \text{otherwise.} \end{cases}$$

Then $J$ is convex, proper and lower semicontinuous (see Barbu [1]).

Under the above assertions we can reformulate the problem $(P')$ in terms of variational inequalities as follows: find $u \in W_0^{1,p}(\Omega)$ such that

$$\int_\Omega a(x, \nabla u) \cdot \nabla(w - u)\,dx + J(w) - J(u) \geq \int_\Omega G(u)(w - u)\,dx \qquad \text{for all } w \in W_0^{1,p}(\Omega)\,. \tag{5}$$

LEMMA 2. *Let hypotheses $(H_1)$ and $(H_2)$ be fulfilled. Then $u \in [\underline{u}, \overline{u}]$ is a solution of (5) if and only if $u$ is a solution of the problem $(P')$.*

*Proof.* Let $u \in [\underline{u}, \overline{u}]$ satisfy the variational inequality (5). Then

$$J(w) \geq J(u) + \int_\Omega G(u) \cdot (w - u)\,dx - \int_\Omega a(x, \nabla u) \cdot \nabla(w - u)\,dx.$$

It follows that

$$\operatorname{div}(a(x, \nabla u)) + G(u) \in \partial J(u) \qquad \text{in } W^{-1,p'}(\Omega). \tag{6}$$

It follows by Brezis [2] [Corollaire 1] that any subgradient $v \in \partial J(u)$ of the functional $J : W_0^{1,p}(\Omega) \to (-\infty, +\infty]$ at $u \in W_0^{1,p}(\Omega)$ belongs to $L^1(\Omega)$ and satisfies

$$v(x) \in \partial j(u(x)) = \beta(u(x)) \qquad \text{a.e. in } \Omega\,. \tag{7}$$

Furthermore

$$h^-(\underline{u}(x)) \leq h^-(u(x)) \leq \beta(u(x)) \leq h^+(\overline{u}(x)) \leq h^+(\overline{u}(x)) \quad \text{a.e. in } \Omega.$$

Thus

$$|v| \leq |H^+(\overline{u})| + |H^-(\underline{u})|\,. \tag{8}$$

By $(H_2)$, the right-hand side of (8) belongs to $L^{p'}(\Omega)$. It follows that $v \in L^{p'}(\Omega)$. Thus there exists $v \in L^{p'}(\Omega)$ such that

$$\operatorname{div}(a(x, \nabla u)) + G(u) = v \qquad \text{in } W^{-1,p'}(\Omega)$$

or, equivalently,

$$\int_\Omega a(x, \nabla u) \cdot \nabla w \, dx + \int_\Omega v \cdot w \, dx = \int_\Omega G(u) w \, dx \quad \text{for all } w \in W_0^{1,p}(\Omega). \tag{9}$$

Relations (7) and (9) imply that $u \in W_0^{1,p}(\Omega)$ is a solution of the problem $(P')$.

Conversely, let $u \in [\underline{u}, \overline{u}]$ be a solution of the problem $(P')$. Then there exists $v \in L^{p'}(\Omega)$ such that $v \in \beta(u) = \partial j(u(x))$ and the relation (9) is fulfilled. Since $v(x) \in \partial j(u(x))$ we have

$$j(s) \geq j(u(x)) + v(x)(s - u(x)) \, . \tag{10}$$

Taking $s = 0$ in (10) we obtain, by means of (3) and (4) that $0 \leq j(u(x)) \leq v(x)u(x)$. Thus

$$j(u(\cdot)) \in L^1(\Omega) \quad \text{and} \quad J(u) = \int_\Omega j(u(x)) \, dx \, . \tag{11}$$

Let $w \in W_0^{1,p}(\Omega)$. Taking $s = w(x)$ in (10) we obtain

$$\int_\Omega j(w(x)) \, dx - \int_\Omega j(u(x)) \, dx \geq \int_\Omega v(x)(w(x) - u(x)) \, dx \, . \tag{12}$$

From (9), substituting $w$ by $w - u \in W_0^{1,p}(\Omega)$ we get, by means of (12)

$$\int_\Omega a(x, \nabla u) \cdot \nabla (w - u) \, dx + J(w) - J(u) \geq \int_\Omega G(u) \cdot (w - u) \, dx \quad \text{for all } w \in W_0^{1,p}(\Omega).$$

This means that $u$ is a solution of the variational inequality (5). □

*Remark* 1. If $u$ is a solution of $(P')$ then, by (11), $J(u) < +\infty$. The result also holds also if we replace $u$ by a super-solution $\overline{u}$ or by a sub-solution $\underline{u}$.

Set $v^+ = \max\{v, 0\}$.

LEMMA 3. *Let $u$, $v \in L^p(\Omega)$ such that $J(u)$ and $J(v)$ are finite. Then*

$$J(u - (u - v)^+) - J(u) + J(v + (u - v)^+) - J(v) = 0. \tag{13}$$

*Proof.* Let $\Omega_+ := \{x \in \Omega |\, u > v\}$ and $\Omega_- := \{x \in \Omega |\, u \le v\}$. Since $(u-v)^+ = 0$ in $\Omega_-$ and $(u-v)^+ = u - v$ in $\Omega_+$ we obtain

$$J(u-(u-v)^+) = \int_{\Omega_+} j(v)\,dx + \int_{\Omega_-} j(u)\,dx \le \infty \tag{14}$$

$$J(v+(u-v)^+) = \int_{\Omega_+} j(u)\,dx + \int_{\Omega_-} j(v)\,dx \le \infty \tag{15}$$

By (14) and (15) we obtain (13).

Consider now the following variational inequality: given $z \in L^p(\Omega)$, find $u \in W_0^{1,p}(\Omega)$ such that

$$\int_\Omega a(x,\nabla u)\cdot\nabla(w-u) + J(w) - J(u) \ge \int_\Omega G(z)(w-u)\,dx \tag{16}$$

$$\text{for all } w \in W_0^{1,p}(\Omega).$$

The variational inequality (16) defines a mapping $T : z \to u$ and each fixed point of $T$ yields a solution of (5) and conversely.

LEMMA 4. *Let hypotheses $(H_1)$ and $(H_2)$ be satisfied. Then for each $z \in [\underline{u},\overline{u}]$ the variational inequality (16) has a unique solution $u = Tz \in [\underline{u},\overline{u}]$. Moreover, there is a constant $C > 0$ such that $\|Tz\|_{W_0^{1,p}(\Omega)} \le C$, for any $z \in [\underline{u},\overline{u}]$.*

*Proof. Existence.* Let $z \in [\underline{u},\overline{u}]$ be arbitrarily given. Then $G(z)$ is measurable and $G(z) \in L^{p'}(\Omega)$, due to the estimate

$$|G(z)| \le |G(\overline{u})| + |G(\underline{u})|$$

and after observing that the right-hand side of the above inequality is in $L^{p'}(\Omega)$, by $(H_2)$.

We now apply Theorem II.8.5 in Lions [5]. We first observe that the above assertions show that the mapping $W_0^{1,p}(\Omega) \ni u \to \int_\Omega G(z)u$ is in $W^{-1,p'}(\Omega)$.

Consider the Leray-Lions operator $A : W_0^{1,p}(\Omega) \to W^{-1,p'}(\Omega)$ defined by

$$\langle Au, w\rangle = \int_\Omega a(x,\nabla u)\cdot\nabla w\,dx\,.$$

We show that $A$ is a pseudo-monotone operator. For this aim it is enough to prove that $A$ is bounded, monotone and hemi-continuous (see Lions [5] [Prop. II.2.5]).

Condition $(a_3)$ yields the boundedness of $A$. Indeed

$$\|Au\|_{W^{-1,p'}}(\Omega) \leq C(\|k\|_{L^{p'}(\Omega)} + \|\nabla u\|_{L^p(\Omega)}^{p-1}) .$$

We also observe that $(a_2)$ implies that $A$ is a monotone operator.

In order to justify the hemi-continuity of $A$, let us consider a sequence $(\lambda_n)_{n\geq 1}$ converging to $\lambda$. Then, for given $u,\ v,\ w\ \in W_0^{1,p}(\Omega)$, we have

$$a(x, \nabla(u + \lambda_n v)) \cdot \nabla w \to a(x, \nabla(u + \lambda v)) \cdot \nabla w \qquad \text{a.e. in } \Omega .$$

From the boundedness of $\{\lambda_n\}$ and condition $(a_3)$ we obtain that the sequence $\{|a(x, \nabla(u + \lambda_n v))\nabla w|\}$ is bounded by a function which belongs to $L^1(\Omega)$. Using the Lebesgue dominated convergence theorem it follows that

$$\langle A(u + \lambda_n v), w\rangle \to \langle A(u + \lambda_n v), w\rangle \qquad \text{as } n \to \infty.$$

Hence the application $\lambda \to \langle A(u + \lambda v,\ w\rangle)$ is continuous.

It follows that all assumptions of Theorem II.8.5 in [5] are fulfilled, so the problem (16) has at least a solution.

*Uniqueness.* Let $u_1$ and $u_2$ be two solutions of (16). Then taking $w = u_2$ as a test function for the solution $u_1$, we obtain

$$\int_\Omega a(x, \nabla u_1) \cdot \nabla(u_2 - u_1)\, dx + J(u_2) - J(u_1) \geq \int_\Omega G(z)(u_2 - u_1)\, dx .$$

Similarly we find

$$\int_\Omega a(x, \nabla u_2) \cdot \nabla(u_1 - u_2)\, dx + J(u_1) - J(u_2) \geq \int_\Omega G(z)(u_1 - u_2)\, dx .$$

Therefore

$$\int_\Omega (a(x, \nabla u_1) - a(x, \nabla u_2)) \cdot (\nabla u_1 - \nabla u_2)\, dx \leq 0 .$$

So, by $(a_2)$, it follows that $\nabla u_1 = \nabla u_2$, so $u_1 = u_2 + C$ in $\Omega$. Since $u_1 = u_2 = 0$ on $\partial\Omega$, it follows that $u_1 = u_2$ in $\Omega$.

From (3) and (4) we deduce that $J(0) = 0$ and $J(u) \geq 0$. Moreover, the variational inequality (16) implies

$$\int_\Omega a(x, \nabla u) \cdot \nabla(-u)\, dx + J(0) - J(u) \geq -\int_\Omega G(z) u\, dx .$$

Thus

$$\int_\Omega a(x, \nabla u) \cdot \nabla u \, dx \le \int_\Omega G(z) u \, dx \, .$$

This last inequality, assumption $(a_1)$ and Hölder's inequality yield

$$\begin{aligned}\lambda \cdot \|u\|^p_{W_0^{1,p}(\Omega)} &\le \int_\Omega G(z) u \, dx \le \|G(z)\|_{L^{p'}(\Omega)} \cdot \|u\|_{L^p} \\ &\le C_1 \left(\|G(\overline{u})\|_{L^{p'}(\Omega)} + \|G(\underline{u})\|_{L^{p'}(\Omega)}\right) \|u\|_{W_0^{1,p}(\Omega)} \, .\end{aligned}$$

Thus $u = Tz$ verifies

$$\|u\|^{p-1}_{W_0^{1,p}(\Omega)} \le C_1(\|G(\overline{u})\|_{L^{p'}(\Omega)} + \|G(\underline{u})\|_{L^{p'}(\Omega)}) = C_2 \, .$$

This implies that there exists a universal constant $C > 0$ such that

$$\|u\|_{W_0^{1,p}(\Omega)} \le C \, .$$

So, in order to conclude our proof, it is enough to show that $u \in [\underline{u}, \overline{u}]$. But, by the definition of an upper solution, there exists $\overline{v} \in L^{p'}(\Omega)$ such that $\overline{v} \in \beta(\overline{u}(x))$ and

$$\int_\Omega a(x, \nabla \overline{u}) \cdot \nabla w \, dx + \int_\Omega \overline{v} \cdot w \, dx \ge \int_\Omega G(\overline{u}) w \, dx, \tag{17}$$

$$\text{for all } w \in W_0^{1,p}(\Omega) \cap L^p_+(\Omega).$$

The solution $u = Tz$ of the variational inequality (16) satisfies

$$\int_\Omega a(x, \nabla u) \cdot \nabla (w - u) \, dx + J(w) - J(u) \ge \int_\Omega G(z)(w - u) \, dx \tag{18}$$

$$\text{for all } w \in W_0^{1,p}(\Omega).$$

Setting $\overline{v} \in \beta(\overline{u}) = \partial j(\overline{u})$, we have

$$j(s) \ge j(\overline{u}(x)) + \overline{v}(x)(s - \overline{u}(x)) \quad \text{for all } \ s \in \mathbf{R}. \tag{19}$$

Taking $s := \overline{u}(x) + (u(x) - \overline{u}(x))^+$ in (19) we find by integration

$$J(\overline{u} + (u - \overline{u})^+) \ge J(\overline{u}) + \int_\Omega \overline{v}(u - \overline{u})^+ \, dx \, . \tag{20}$$

Choosing now $w = (u - \overline{u})^+$ in (17) we obtain

$$\int_\Omega a(x, \nabla\overline{u}) \cdot \nabla(u - \overline{u})^+ \, dx + \int_\Omega \overline{v} \cdot (u - \overline{u})^+ \, dx \geq \int_\Omega G(\overline{u}) \cdot (u - \overline{u})^+ \, dx \,. \tag{21}$$

Relations (20) and (21) yield

$$\int_\Omega a(x, \nabla\overline{u}) \nabla(u - \overline{u})^+ \, dx + J(\overline{u} + (u - \overline{u})^+) - J(\overline{u}) \geq \int_\Omega G(\overline{u}) \cdot (u - \overline{u})^+ \, dx. \tag{22}$$

Taking $w = u - (u - \overline{u})^+$ in (18), we obtain

$$\int_\Omega a(x, \nabla u) \cdot (-\nabla(u - \overline{u})^+) \, dx + J(u - (u - \overline{u})^+) - J(u) \geq - \int_\Omega G(z)(u - \overline{u})^+ \, dx.$$

Since $z \in [\underline{u}, \overline{u}]$ and $G : L^p(\Omega) \to L^p(\Omega)$ is nondecreasing, it follows that

$$\begin{aligned} \int_\Omega a(x, \nabla u) \cdot \nabla(u - \overline{u})^+ dx + J(u - (u - \overline{u})^+) - J(u) \\ \geq - \int_\Omega G(\overline{u})(u - \overline{u})^+ \, dx. \end{aligned} \tag{23}$$

From (22), (23) and Lemma 3 we have

$$\int_\Omega (a(x, \nabla u) - a(x, \nabla\overline{u}) \cdot \nabla(u - \overline{u})^+ \, dx \leq 0. \tag{24}$$

Let $\Omega_+ = \{x \in \Omega \,|\, u \leq \overline{u}\}$ and $\Omega_- = \{x \in \Omega \,|\, u > \overline{u}\}$. Since $(u - \overline{u})^+ = 0$ in $\Omega_+$ and $(u - \overline{u})^+ = u - \overline{u}$ in $\Omega_-$, it follows by (24) that

$$\int_{\Omega_-} (a(x, \nabla u) - a(x, \nabla\overline{u}) \cdot \nabla(u - \overline{u})^+ \, dx \leq 0 \,.$$

So, by $(a_2)$ and the definition of $\Omega_-$, we obtain meas $(\Omega^-) = 0$, hence $u \leq \overline{u}$ a.e. in $\Omega$. Proceeding in the same way we prove that $\underline{u} \leq u$. □

LEMMA 5. *The operator $T$ defines a monotone nondecreasing mapping from $[\underline{u}, \overline{u}]$ to $[\underline{u}, \overline{u}]$.*

*Proof.* Let $z_1,\ z_2 \in [\underline{u}, \overline{u}]$ be such that $z_1 \le z_2$. By Lemma 4, we obtain that $Tz_1,\ Tz_2 \in [\underline{u}, \overline{u}]$ and

$$\int_\Omega a(x, \nabla Tz_1) \cdot \nabla(w - Tz_1)\, dx + J(w) - J(Tz_1) \geq \int_\Omega G(z_1)(w - Tz_1)\, dx \tag{25}$$

$$\int_\Omega a(x, \nabla Tz_2) \cdot \nabla(w - Tz_2)\, dx + J(w) - J(Tz_1) \geq \int_\Omega G(z_2)(w - Tz_2)\, dx\,. \tag{26}$$

Taking $w = Tz_1 - (Tz_1 - Tz_2)^+$ in (25) and $w = Tz_2 + (Tz_1 - Tz_2)^+$ in (26), we get

$$-\int_\Omega a(x, \nabla Tz_1) \cdot \nabla(Tz_1 - Tz_2)^+\, dx + J(Tz_1 - (Tz_1 - Tz_2)^+) - J(Tz_1) \geq \int_\Omega G(z_1)(-(Tz_1 - Tz_2)^+)\, dx$$

$$\int_\Omega a(x, \nabla Tz_2) \cdot \nabla(Tz_1 - Tz_2)^+\, dx + J(Tz_2 + (Tz_1 - Tz_2)^+) - J(Tz_2) \geq \int_\Omega G(z_2)(Tz_1 - Tz_2)^+\, dx\,.$$

Summing up these inequalities we get, by means of (13),

$$\int_\Omega (a(x, \nabla Tz_1) - a(x, \nabla Tz_2)) \cdot \nabla(Tz_1 - Tz_2)^+\, dx \leq \int_\Omega (G(z_1) - G(z_2))(Tz_1 - Tz_2)^+\, dx\,.$$

But $G(z_1) \le G(z_2)$, since $G$ is a nondecreasing operator. Therefore, by the above inequality we obtain

$$\int_\Omega a(x, (\nabla Tz_1) - a(x, \nabla Tz_2)) \cdot \nabla(Tz_1 - Tz_2)^+ \, dx \le 0 \, .$$

With the same argument as for proving (24) we obtain $Tz_1 \le Tz_2$. □

*Proof of Theorem 1 completed.* Assume that $g$ is right continuous. Define

$$u^{n+1} = Tu^n \, , \tag{27}$$

where $u^0 = \overline{u}$. Then, by Lemma 4, $\{u^n\}$ is nondecreasing, $u^n \in [\underline{u}, \overline{u}]$, and there is a constant $C$ such that

$$\|u^n\|_{W_0^{1,p}(\Omega)} \le C. \tag{28}$$

The compact embedding $W_0^{1,p}(\Omega) \hookrightarrow L^p(\Omega)$ and (28) ensure that there exists $u \in W_0^{1,p}(\Omega)$ such that, up to a subsequence,

$$\begin{aligned} u^n &\to u \quad \text{strongly in } L^p(\Omega) \\ u^n &\rightharpoonup u \quad \text{weakly in } W_0^{1,p}(\Omega) \\ u_n &\to u \quad \text{a.e. in } \Omega. \end{aligned}$$

By Lemma 4, there exists $u' \in W_0^{1,p}(\Omega)$, $u' \in [\underline{u}, \overline{u}]$ such that $u' = Tu$. We prove in what follows that $u$ is a fixed point of $T$ i.e. $u' = u$.

From (27) and by the definition of $T$ we obtain

$$\int_\Omega a(x, \nabla u^{n+1}) \nabla(w - u^{n+1}) dx + J(w) - J(u^{n+1}) \ge \int_\Omega G(u^n)(w - u^{n+1}) \tag{29}$$

for all $w \in W_0^{1,p}(\Omega)$.

Also, from $Tu = u'$, we have

$$\int_\Omega a(x, \nabla u') \nabla(w - u') \, dx + J(w) - J(u') \ge \int_\Omega G(u) \cdot (w - u') \, dx \tag{30}$$

for all $w \in W_0^{1,p}(\Omega)$.

Taking $w = u'$ in (29) and $w = u^{n+1}$ in (30), we get

$$\int_\Omega a(x, \nabla u^{n+1}) \nabla(u' - u^{n+1}) \, dx + J(u') - J(u^{n+1}) \ge \int_\Omega G(u^n) \cdot (u' - u^{n+1}) \, dx$$

$$\int_\Omega a(x,\nabla u')\nabla(u^{n+1}-u')\,dx + J(u^{n+1}) - J(u') \geq \int_\Omega G(u)\cdot(u^{n+1}-u')\,dx\,.$$

So, by (29) and (30), $J(u') < \infty$ and $J(u^{n+1}) < \infty$. Summing up the last two inequalities we obtain

$$\int_\Omega (a(x,\nabla u') - a(x,\nabla u^{n+1})\cdot\nabla(u'-u^{n+1})\,dx \leq \int_\Omega (G(u)-G(u^n))\,(u'-u^{n+1})\,dx. \tag{31}$$

Since $G$ is right continuous we have $G(u^n) \to G(u)$ in $\Omega$. We also have

$$|G(u)-G(u^n)|\left(u-u^{n+1}\right) \leq 2\left(|G(\underline{u})|+|G(\overline{u})|\right)\left(|\underline{u}|+|\overline{u}|\right) \in L^1(\Omega).$$

By $(a_2)$ and the Lebesgue dominated convergence theorem, we deduce from (31) that

$$\int_\Omega (a(x,\nabla u') - a(x,\nabla u^n)\cdot\nabla(u'-u^n)\,dx \to 0\,. \tag{32}$$

This implies that $\nabla u^n \to \nabla u'$ a.e. in $\Omega$.

Relation (32) implies that (up to a subsequence)

$$(a(x,\nabla u') - a(x,\nabla u^n))\cdot\nabla(u'-u^n) \to 0 \qquad \text{a.e. } x \in \Omega. \tag{33}$$

This leads to $\nabla u^n \to \nabla u'$ a.e. in $\Omega$. Indeed, if not, there exists $x \in \Omega$ such that (up to a subsequence), $\nabla u^n(x) \to \xi \in \overline{\mathbf{R}}^N$ for $\xi \neq \nabla u'$. Passing to the limit in (33) we obtain

$$(a(x,\nabla u') - a(x,\xi))\cdot(\nabla u' - \xi) = 0\,,$$

which contradicts $(a_2)$. So, we have proved that $\nabla u^n \to \nabla u$. Using the fact that $u^n \rightharpoonup u$ weakly in $W_0^{1,p}(\Omega)$, we conclude that $\nabla u' = \nabla u$, thus $u' = u$. Replacing $u'$ by $u$ in (30) we get

$$\int_\Omega a(x,\nabla u)\cdot\nabla(w-u)\,dx + J(w) - J(u) \geq \int_\Omega G(u)(w-u)\,dx \quad \text{for all } w \in W_0^{1,p}(\Omega).$$

Hence $u$ is a fixed point of $T$ and a solution for the problem $(P')$.

In order to prove that $u$ is a maximal solution of (3) with respect to the order interval $[\underline{u},\overline{u}]$, take any other solution $\hat{u} \in [\underline{u},\overline{u}]$ of the problem $(P')$.

Then $\hat{u}$ is in particular a sub-solution satisfying $\hat{u} \le \bar{u}$. Starting again the iteration (27) with $u^0 = \bar{u}$ we obtain

$$\hat{u} \le \cdots \le u^{n+1} \le u^n \le \cdot \le u^0 = \bar{u} .$$

It follows that $\hat{u} \le u$, which concludes our proof. □

Simona Dăbuleanu – Vicenţiu Rădulescu
*Department of Mathematics*
*University of Craiova*
*1100 Craiova, Romania*
*E-mail: vicrad@yahoo.com*